\documentclass{amsart}
\usepackage{amsthm}                 %Definitionen, Theoreme etc.
\usepackage{amsmath}                %siehe userguide, z.B. boxed{}, text-Befehl etc.
\usepackage{amsfonts,amssymb}       %ein paar zusätzliche Symbole
\usepackage{graphicx}               %Bilder
\usepackage{float}                  %Bilder/Positionierung
\usepackage{rotating}               %Bilder/Positionierung
\restylefloat{figure}               %Bilder/Positionierung
\usepackage{mathrsfs}               %mathscr-Buchstaben
\usepackage{fancyhdr}               %Kopfzeile
\usepackage{hyperref}               %Referenzen anklicken
\usepackage[numbers,sort&compress]{natbib}
\usepackage{color}
\theoremstyle{plain}
\newtheorem{lem}{Lemma}
\newtheorem{prop}[lem]{Proposition}
\newtheorem{thm}[lem]{Theorem}

\theoremstyle{definition}

\newcommand{\R}{\mathbb R}

\newcommand{\Z}{\mathbb Z}
\newcommand{\N}{\mathbb N}

\renewcommand{\L}{\mathcal L}

\renewcommand{\d}{\,\text{\rm d}}
\newcommand{\dw}{\text{\rm d}}

\renewcommand{\S}{\mathbb S}
\renewcommand{\phi}{\varphi}

\newcommand{\eps}{\varepsilon}

\newcommand{\set}[2]{\left\{#1;\;#2\right\}}

\newcommand{\bea}{\begin{eqnarray}}
\newcommand{\eea}{\end{eqnarray}}
\newcommand{\beq}{\begin{equation}}
\newcommand{\eeq}{\end{equation}}
\renewcommand{\phi}{\varphi}

\renewcommand{\autoref}[1]{\text{Eq.}~\eqref{#1}}

\begin{document}
\title{Analysis of a mathematical model for the growth of cancer cells}
\author{Martin Kohlmann}
\address{Peter L. Reichertz Institute for Medical Informatics, University of Braunschweig, D-38106 Braunschweig, Germany}
\email{martin.kohlmann@plri.de}
\begin{abstract}
In this paper, a two-dimensional model for the growth of multi-layer tumors is presented. The model consists of a free boundary problem for the tumor cell membrane and the tumor is supposed to grow or shrink due to cell proliferation or cell dead. The growth process is caused by a diffusing nutrient concentration $\sigma$ and is controlled by an internal cell pressure $p$. We assume that the tumor occupies a strip-like domain with a fixed boundary at $y=0$ and a free boundary $y=\rho(x)$, where $\rho$ is a $2\pi$-periodic function. First, we prove the existence of solutions $(\sigma,p,\rho)$ on a scale of small H\"older spaces and show that our model allows for flat stationary solutions. As a main result we establish that these equilibrium points are locally asymptotically stable under small perturbations.
\end{abstract}
\maketitle
\tableofcontents
\section{Introduction}\label{sec_intro}
The mathematical modeling of cancer growth is a challenging area of research in the applied sciences nowadays. The complex growth process of a tumor cell can be captured within different mathematical models, e.g., models consisting of a system of coupled partial differential equations which arise from reaction-diffusion equations and a mass conservation law \cite{BC95,BC96,G72,WK97,WK98}. In this context, tumor growth is often considered as a free boundary problem \cite{BF03,BF03',BC97,CCF05,C07,CuiE09,CF02,CF03,CW05,Xu11}. Sometimes, it is also justified to treat cancer growth as an incompressible flow in a porous medium so that cells move in accordance with Darcy's law, see, e.g., \cite{CCF05,CuiE09,F05}. We refer the reader to the review papers \cite{AP09, BFG06, Letal10} which present a variety of other tumor growth models.

Over the years, the following general non-dimensionalized moving boundary value problem
\beq
\left\{
\begin{array}{rclc}
  \Delta\sigma     & = & f(\sigma)                   & \text{in }\Omega(t),          \\
  \Delta p         & = & -g(\sigma)                  & \text{in }\Omega(t),          \\
  \sigma           & = & \bar\sigma                  & \text{on }\partial\Omega(t),  \\
  p                & = & \gamma\kappa                & \text{on }\partial\Omega(t),  \\
  V                & = & -\partial_\nu p             & \text{on }\partial\Omega(t),
\end{array}
\right.
\label{introduction_problem}\eeq
has turned out to be very appropriate to describe the growth of tumor cells; see, e.g., \cite{F05}. Here $\sigma=\sigma(t,x)$ and $p=p(t,x)$ are unknown functions defined on the time-space manifold $\bigcup_{t\geq 0}\{t\}\times\Omega(t)$ with an a priori unknown time-dependent domain $\Omega(t)\subset\R^n$. The boundary $\partial\Omega(t)$ has to be determined together with the functions $\sigma$ and $p$ which model the nutrient concentration and the internal pressure for the tumor cell described by $\Omega(t)$. Furthermore, $\bar\sigma>0$ is a constant and $\kappa$ and $\partial_\nu$ denote the mean curvature and the normal derivative (with respect to the outward normal) for the boundary $\partial\Omega$. Finally, $\gamma>0$ is the surface tension coefficient and $V$ denotes the normal velocity of the free boundary. In \eqref{introduction_problem}, the tumor $\Omega(t)$ receives a constant supply of nutrient on its boundary and the pressure on $\partial\Omega(t)$ is proportional to the curvature of the free surface $\partial\Omega(t)$, as proposed in \cite{G76}. The evolution equation for the moving boundary comes from an application of Darcy's law from fluid mechanics. We also impose the initial condition $\Omega(0)=\Omega_0$, where $\Omega_0\subset\R^n$ is a given bounded domain in $\R^n$ with a sufficiently smooth boundary.

Typical choices for $f$ and $g$ are
$$f(\sigma)=\lambda\sigma,\quad g(\sigma)=\mu(\sigma-\tilde\sigma)$$
where $\lambda,\mu,\tilde\sigma>0$ are constants, \cite{BC95}, or
\bea f(\sigma) & = & \beta\frac{A\sigma^{m_1}}{\sigma_c^{m_1}+\sigma^{m_1}}+h(\sigma),\nonumber\\
g(\sigma) & = & \frac{A\sigma^{m_1}}{\sigma_c^{m_1}+\sigma^{m_1}}-B\left(1-\frac{\delta\sigma^{m_2}}{\sigma_d^{m_2}+\sigma^{m_2}}\right),
\nonumber\eea
with positive constants $A$, $B$, $\beta$, $\delta$, $m_1$, $m_2$, $\sigma_c$, $\sigma_d$ and a non-negative increasing function $h$, \cite{WK97}. In \cite{CE07}, the authors assume that $f$ and $g$ are general functions satisfying
\begin{itemize}
\item $f,g\in C^{\infty}[0,\infty)$,
\item $f'(\sigma)>0$ for $\sigma\geq 0$ and $f(0)=0$,
\item $g'(\sigma)>0$ for $\sigma\geq 0$ and $g(\tilde\sigma)=0$ for some $\tilde\sigma>0$,
\item $\tilde\sigma<\bar\sigma$.
\end{itemize}
In this paper, we will suppose that
$$f(\sigma)=\sigma,\quad g(\sigma)=\mu(\sigma-\tilde\sigma),$$
with positive parameters $\mu,\tilde\sigma$. Since $\Delta p$ is minus the divergence of the cell velocity field, the meaning of $\Delta p=-\mu(\sigma-\tilde\sigma)$ is that tumor volume is produced, if $\sigma$ is above the proliferation threshold $\tilde\sigma$, and that the tumor volume decreases in the opposite case.

The tumor growth model which we study in this paper has the following form: The tumor is assumed to occupy a two-dimensional region of the form
$$\Omega_\rho(t):=\set{(x,y)\in\R^2}{0<y<\rho(t,x)}$$
where $t$ is the time variable and $\rho(t,x)$ is an unknown positive $2\pi$-periodic function. The upper boundary of the tumor is denoted as
$$\Gamma_\rho(t):=\set{(x,y)\in\R^2}{y=\rho(t,x)},$$
its lower boundary is $\Gamma_0=\set{(x,y)\in\R^2}{y=0}$. A similar situation is discussed in \cite{CuiE09} where the authors explain that the strip-shaped model refers to the growth of multi-layer tumors, a kind of in vitro tumors cultivated in laboratory by using the recently developed tissue culture technique, \cite{KSoH04,KCM99,M97}. While there are only a few works dealing with strip-shaped domains, a variety of papers considering radially symmetric models for tumor growth have been published, cf., e.g., the seminal paper \cite{FR99}.
\begin{figure}[H]
\begin{center}
\includegraphics[width=8cm]{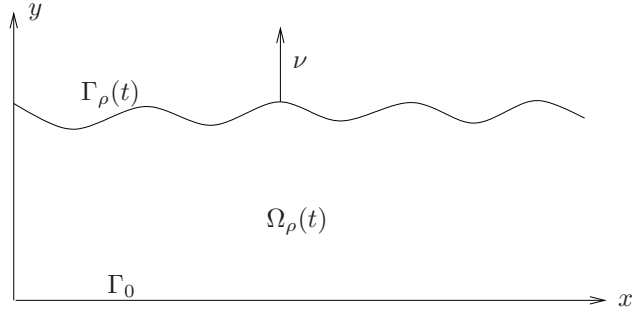}
\caption{A free boundary problem modeling multi-layer tumors.}
\end{center}
\end{figure}
\begin{picture}(0,0)
\put(160,140){$\nu$}
\put(150,80){$\Omega_\rho(t)$}
\put(90,55){$\Gamma_0$}
\put(80,128){$\Gamma_\rho(t)$}
\put(283,50){$x$}
\put(60,160){$y$}
\end{picture}
The following conditions on the tumor growth process are imposed: We assume that the tumor is constantly supplied with nutrient on $\partial\Omega_\rho(t)$; precisely, the nutrient concentration is $\bar\sigma_2>0$ on $\Gamma_\rho(t)$ and $\bar\sigma_1>0$ on $\Gamma_0$. As in \eqref{introduction_problem}, we assume that the pressure on $\Gamma_\rho$ is given by $\gamma\kappa_{\Gamma_\rho}$. Since tumor cells should only grow in the positive $y$ direction, we assume that $p_y=0$ on $\Gamma_0$. Thus we are led to study the following system of equations:
\beq
\left\{
\begin{array}{rclcc}
  \Delta\sigma     & = & \sigma                      & \text{in }\Omega_\rho(t), & t>0, \\
  \Delta p         & = & -\mu(\sigma-\tilde\sigma)                      & \text{in }\Omega_\rho(t), & t>0, \\
  \rho_t           & = & - \frac{\partial p}{\partial\nu}       & \text{on }\Gamma_\rho(t), & t>0, \\
  \sigma           & = & \bar\sigma_{2}             & \text{on }\Gamma_\rho(t), & t>0, \\
  \sigma           & = & \bar\sigma_{1}                           & \text{on }\Gamma_0,       & t>0, \\
  p                & = & \gamma\kappa_{\Gamma_\rho}  & \text{on }\Gamma_\rho(t), & t>0, \\
  p_y              & = & 0                        & \text{on }\Gamma_0,       & t>0, \\
  \rho             & = & \rho_0                     & \text{for } t=0, &
\end{array}
\right.
\label{problem}\eeq
where $\nu=(-\rho_x,1)$ denotes the outward normal on $\Gamma_\rho(t)$ with respect to $\Omega_\rho(t)$ and $\rho_0$ is a given periodic function. The third equation in \eqref{problem} can be derived as follows: we assume that the normal velocity $V$ of the boundary $\Gamma_\rho$ is equal to the cell movement velocity in the direction of the outward unit normal
$$\nu_0=\frac{1}{\sqrt{1+\rho_x^2}}(-\rho_x,1),$$
cf. \cite{BF03,BC95,E04,G72,WK97}. Then Darcy's law and the relation $V=\rho_t(1+\rho_x^2)^{-1/2}$ (cf.~\cite{E04}) imply that $V=-\nu_0\cdot\nabla p$ and that the motion of the free surface $\Gamma_\rho$ is modeled by the third equation of \eqref{problem}.

In the model presented in \cite{CuiE09}, the lower boundary $\Gamma_0$ is supposed to be impermeable for glucose and oxygen. The authors also comment briefly on a variant of their model obtained by exchanging the boundary conditions for $\sigma$ on $\Gamma_\rho$ and $\Gamma_0$. The novel aspect of the problem presented in the paper at hand is that we allow for supply of nutrient on both boundary components.

% Related mathematical models have been discussed in \cite{CuiE09,MK11,ZECui08,ZECui08b}. Our approach gets in line with a wide rage of papers discussing mathematical models for cancer growth; many of them discussing spherically symmetric tumors as in the seminal papers \cite{BC95,FR99,WK97}. In addition, the model under discussion is strongly related to the Hele-Shaw problem (as explained in the introduction of \cite{F05} and in \cite{ES97}), i.e., the growth of tumor cells is captured by a similar approach as the movement of a viscous fluid.

Our paper can be outlined as follows:
In Section~\ref{sec_wdh}, we recall some elementary facts, definitions and notation from \cite{CuiE09} which will be important for our approach to the system \eqref{problem}. In Section~\ref{sec_sol}, we prove that the system~\eqref{problem} is well-posed on a scale of small H\"older spaces and we compute its flat stationary solutions. Finally, in Section~\ref{sec_stabil}, we linearize the system~\eqref{problem} at such an equilibrium and prove that it is asymptotically stable under small perturbations.
\\[.1cm]

\emph{Acknowledgement.} The author thanks the anonymous referees for helpful suggestions that improved the preliminary version of the paper.
\section{Preliminaries}\label{sec_wdh}
Let $\S=\R/2\pi\Z$. Function spaces of $2\pi$-periodic functions will be identified with the corresponding spaces on $\S$. In the following $C_+([0,T)\times\S)$ stands for the cone of all positive functions in $C([0,T)\times\S)$, for any $T>0$. We will make use of the little H\"older spaces $h^{k+\alpha}(\S)$ which are defined as the closure of $C^{\infty}(\S)$ in the usual H\"older space $C^{k+\alpha}(\S)$, for $k\in\N$ and $\alpha\in(0,1)$. Similarly, for any open set $\Omega\subset\R^2$, we define $h^{k+\alpha}(\overline\Omega)$ as the closure of $C^{\infty}(\overline\Omega)$ in $C^{k+\alpha}(\overline\Omega)$. The little H\"older spaces are Banach algebras under pointwise multiplication and the embedding $h^r(\S)\hookrightarrow h^s(\S)$, for $r>s$, is compact. We shall denote the cone of positive functions in $h^{k+\alpha}(\S)$ by $h^{k+\alpha}_+(\S)$.

Let
$$D_{\rho,T}=\set{(t,x,y)}{t\in[0,T),\,x\in\S,\,0\leq y\leq\rho(t,x)}.$$
We call a triple $(\sigma,p,\rho)$ a solution to the problem \eqref{problem} if
\begin{align}
&(\sigma,p)\in C(D_{\rho,T})\times C(D_{\rho,T}),\nonumber\\
&\sigma(t,\cdot)\in h^{4+\alpha}(\overline{\Omega_\rho(t)}),\,p(t,\cdot)\in h^{2+\alpha}(\overline{\Omega_\rho(t)}),\,t\in(0,T),\nonumber\\
&\rho\in C([0,T),h_+^{3+\alpha}(\S))\cap C^1((0,T),h^{1+\alpha}(\S)),\nonumber\\
&(\sigma,p,\rho)\text{ satisfies \eqref{problem} pointwise on }D_{\rho,T}.\nonumber
\end{align}
A solution $(\sigma,p,\rho)$ is called stationary, if $(\sigma,p,\rho)=(\sigma(x,y),p(x,y),\rho(x))$. We call a stationary solution flat if $(\sigma,p)=(\sigma(y),p(y))$ and $\rho$ is a positive constant. Finally, we call the problem \eqref{problem} well-posed if there exist a time $T>0$ and a unique solution $(\sigma,p,\rho)$ on $[0,T)$ as defined above.

The usual way of treating a free boundary problem like \eqref{problem} is to transform it to a problem on a fixed reference domain. Let $\Omega:=\S\times(0,1)$ with the boundary components $\Gamma_i:=\S\times\{i\}\simeq\S$ for $i\in\{0,1\}$. Given any $\rho\in C^2_+(\S)$, we introduce the map
$$\theta_\rho\colon\overline\Omega\to\overline{\Omega_\rho},\quad(x',y')\mapsto(x',y'\rho(x'));$$
i.e., we will label the coordinates in $\overline{\Omega_\rho}$ by $(x,y)$ and the variables in $\overline\Omega$ by $(x',y')$ with the transformation
$$x=x'\quad\text{and}\quad y=y'\rho(x')\qquad\text{ or }\qquad x'=x\quad\text{and}\quad y'=\frac{y}{\rho(x)}.$$
Since $\rho$ is positive, it follows that $\theta_\rho$ is a $C^2$-diffeomorphism. Let $u\in C(\overline{\Omega_\rho})$ and $v\in C(\overline\Omega)$. Then
$$\theta_\rho^*u:=u\circ\theta_\rho\in C(\overline\Omega)\quad\text{and}\quad(\theta_{\rho})_*v:=v\circ\theta_\rho^{-1}\in C(\overline{\Omega_\rho}).$$
We call $\theta_\rho^*$ the pull-back and $(\theta_\rho)_*$ the push-forward operator for the pair $(\overline\Omega,\overline{\Omega_\rho})$. For $v\in C^2(\overline\Omega)$ and $\rho\in C^2_+(\S)$ we define
$$\mathcal A(\rho)v:=\theta_\rho^*\Delta[(\theta_\rho)_*v],\quad\mathcal B(\rho)v:=\theta_\rho^*(\text{tr}[\nabla(\theta_{\rho})_*v]\cdot\nu),$$
where $\text{tr}$ denotes the trace operator with respect to $\Gamma_\rho$. It is easy to derive the explicit formulae
\beq \mathcal A(\rho)v
= v_{x'x'}-v_{x'y'}\frac{2y'\rho_{x'}}{\rho}+v_{y'y'}\frac{1+(y')^2\rho_{x'}^2}{\rho^2}
-v_{y'}y'\frac{\rho_{x'x'}\rho -2\rho_{x'}^2}{\rho^2}\label{operatorA}\eeq
and
\beq\mathcal B(\rho)v=\left(-v_{x'}+v_{y'}\frac{\rho_{x'}}{\rho}\right)\bigg|_{y'=1}\rho_{x'}+\frac{1}{\rho}v_{y'}|_{y'=1}.\label{operatorB}
\eeq
The straightforward calculations leading to \eqref{operatorA} and \eqref{operatorB} are omitted for the convenience of the reader. We conclude that $\mathcal A$ is uniformly elliptic in $\overline\Omega$, as defined in \cite{GT77}. Using the fact that the little H\"older spaces are Banach algebras, it follows that
$$\mathcal A\in C^{\infty}(h_+^{3+\alpha}(\S);\L(h^{2+i+\alpha}(\overline\Omega),h^{i+\alpha}(\overline\Omega))),\quad i\in\{0,1\}.$$
and
$$\mathcal B\in C^{\infty}(h_+^{3+\alpha}(\S);\L(h^{2+i+\alpha}(\overline\Omega),h^{1+i+\alpha}(\S))),\quad i\in\{0,1\}.$$
In terms of the variables
$$\tau(t):=\theta_{\rho(t)}^*\sigma(t,\cdot)\quad\text{and}\quad q(t):=\theta_{\rho(t)}^*p(t,\cdot),$$
the system \eqref{problem} is equivalent to
\beq
\left\{
\begin{array}{rclcc}
  \mathcal A(\rho(t))\tau      & = & \tau                           & \text{in }\Omega,   & t>0,\\
  \mathcal A(\rho(t))q         & = & -\mu(\tau-\tilde\sigma)        & \text{in }\Omega,   & t>0,\\
  \rho_t                       & = & -\mathcal B(\rho(t))q          & \text{on }\Gamma_1, & t>0,\\
  \tau                         & = & \bar\sigma_2                   & \text{on }\Gamma_1, & t>0,\\
  \tau                         & = & \bar\sigma_1                   & \text{on }\Gamma_0, & t>0,\\
  q                            & = & \gamma\kappa                   & \text{on }\Gamma_1, & t>0,\\
  q_y                          & = & 0                              & \text{on }\Gamma_0, & t>0,\\
  \rho                         & = & \rho_0                         & \text{for }t=0.
\end{array}
\right.
\label{problem2}\eeq
Next, we introduce the following solution operators: The solution of the boundary value problem
$$
\left\{
\begin{array}{rclc}
  \mathcal A(\rho)\tau & = & \tau          & \text{in }\Omega, \\
  \tau                 & = & \bar\sigma_1  & \text{on }\Gamma_0, \\
  \tau                 & = & \bar\sigma_2  & \text{on }\Gamma_1
\end{array}
\right.
$$
is denoted as $\tau=\mathcal R(\rho)(\bar\sigma_1,\bar\sigma_2)$. For given functions $f\in h^{1+\alpha}(\overline\Omega)$ and $k\in h^{2+\alpha}(\S)$, we write the solution of
$$
\qquad
\left\{
\begin{array}{rcll}
  \mathcal A(\rho)q & = & f   & \text{in }\Omega, \\
  q_y               & = & 0   & \text{on }\Gamma_0, \\
  q                 & = & k   & \text{on }\Gamma_1 \\
\end{array}
\right.$$
as $q=\mathcal S(\rho)f+\mathcal T(\rho)k$. We have
$$\mathcal R(\cdot)(\bar\sigma_1,\bar\sigma_2)\in C^{\infty}(h_+^{3+i+\alpha}(\S);h^{3+i+\alpha}(\overline\Omega)),\quad i\in\{0,1\},$$
and
$$(\mathcal S,\mathcal T)\in C^{\infty}(h_+^{3+\alpha}(\S);\mathcal L(h^{1+\alpha}(\overline\Omega),h^{3+\alpha}(\overline\Omega))\times\mathcal L(h^{2+\alpha}(\S),h^{2+\alpha}(\overline\Omega))).$$
Since the sign of $\kappa$ determines the following results in a very significant way, it will be very instructive to remind the reader of the following elementary result: In classical differential geometry one considers a local parametrization $c(x)=(x,\rho(x))$ of the curve $\Gamma_\rho$ and has
$$\kappa(x)=\frac{|c_x(x)\times c_{xx}(x)|}{|c_x(x)|^3}=(1+\rho_x^2)^{-3/2}\rho_{xx}.$$
Note however that the boundary $\Gamma_\rho$ is convex in $(x,\rho(t,x))$ with respect to the outer normal of $\Omega_\rho$ if and only if $\rho(t,x)$ is concave in $x$ and vice versa. This motivates to change the sign of the classical curvature formula so that we will use the identity
$$\kappa=-(1+\rho_x^2)^{-3/2}\rho_{xx}$$
for the curvature of $\Gamma_\rho$. With this modified sign convention we write
$$\kappa_{\Gamma_\rho}=\mathcal P(\rho)\rho,$$
where
$$\mathcal P(\rho)=-(1+\rho_x^2)^{-3/2}\partial_x^2\in C^{\infty}(h_+^{3+\alpha}(\S);\mathcal L(h^{4+\alpha}(\S),h^{2+\alpha}(\S))).$$
\section{Well-posedness and the flat stationary solutions}\label{sec_sol}
In this section, we discuss the question of existence and uniqueness of a solution to \eqref{problem}. In particular, we are interested in stationary solutions to \eqref{problem}, i.e., solutions which do not depend on time. We begin with a proof of the following theorem. In fact, the arguments are just a repetition of what is derived in \cite{CuiE09}.
\thm Given $\rho_0\in h_+^{3+\alpha}(\S)$, there exists $T>0$ such that \eqref{problem} has a solution which is unique in the class $C(D_{\rho,T})\times C(D_{\rho,T})\times C([0,T),h^{3+\alpha}_+(\S))$. Furthermore, if $\rho_0\in h^{4+\alpha}_+(\S)$, then
$\rho\in C([0,T),h_+^{4+\alpha}(\S))\cap C^1([0,T),h^{1+\alpha}(\S))$.
\endthm\rm
\proof Using the notation of Section~\ref{sec_wdh}, it is easy to see that the transformed system~\eqref{problem2} can be written in the form
\beq
\left\{
\begin{array}{rclc}
\rho_t +\Phi(\rho)\rho &=&F(\rho), & t>0,\\
\rho   &=&\rho_0, &t=0,
\end{array}
\right.
\label{abstrevoleqn}\eeq
where
$$\Phi(\rho)=\gamma\mathcal B(\rho)\mathcal T(\rho)\mathcal P(\rho)\quad\text{and}\quad F(\rho)=\mu\mathcal B(\rho)\mathcal S(\rho)\{\mathcal R(\rho)(\bar\sigma_1,\bar\sigma_2)-\tilde\sigma\}.$$
We find that
\beq(\Phi,F)\in C^{\infty}(h_+^{3+\alpha}(\S);\mathcal L(h^{4+\alpha}(\S),h^{1+\alpha}(\S))\times h^{2+\alpha}(\S)).\label{PhiF}\eeq
Since $\gamma>0$, the operator $\Phi(\rho)$ generates, for any $\rho\in h_+^{3+\alpha}(\S)$, a strongly continuous analytic semigroup on the space $h^{1+\alpha}(\S)$ (cf.~Theorem~4.1.~in \cite{ES97}). An application of Amann's local existence, uniqueness and regularity theory for abstract quasilinear evolution equations (see, e.g., Theorem 12.1 and Remarks 12.2 in \cite{A93}) achieves the proof.
\endproof
To find the flat stationary solutions $(\sigma_*(y),p_*(y),\rho_*)$ of \eqref{problem} we first solve the problem
$$
\left\{
\begin{array}{rcl}
  \sigma_*''       & = & \sigma_*, \\
  \sigma_*(0)      & = & \bar\sigma_1, \\
  \sigma_*(\rho_*) & = & \bar\sigma_2,
\end{array}
\right.
$$
and obtain the unique solution
\beq\sigma_*(y)=(\bar\sigma_2-\bar\sigma_1\cosh\rho_*)\frac{\sinh y}{\sinh\rho_*}+\bar\sigma_1\cosh y.\label{sigmastar}\eeq
Next, we observe that the unique solution of
$$
\left\{
\begin{array}{rcl}
  p_*''       & = & -\mu(\sigma_*-\tilde\sigma), \\
  p_*(\rho_*)      & = & 0, \\
  p_*'(0) & = & 0
\end{array}
\right.
$$
is given by
\beq p_*(y)=\mu\frac{\bar\sigma_2-\bar\sigma_1\cosh\rho_*}{\sin\rho_*}(y-\rho_*)
+\mu\left(\bar\sigma_2-\sigma_*(y)-\frac{1}{2}\tilde\sigma(\rho_*^2-y^2)\right).\label{pstar}\eeq
Since we must demand $p'(\rho_*)=0$, we get the condition
\beq\frac{\bar\sigma_1+\bar\sigma_2}{\tilde\sigma}(1-\cosh\rho_*)+\rho_*\sinh\rho_*=0.\label{rhostar}\eeq
Let $\alpha=\frac{\bar\sigma_1+\bar\sigma_2}{\tilde\sigma}$. We suppose
\beq\bar\sigma_1,\bar\sigma_2>\tilde\sigma\label{condbarsigma}\eeq
to obtain that $\alpha>2$. The function
$$f_\alpha\colon(0,\infty)\to\R,\quad x\mapsto\alpha(1-\cosh x)+x\sinh x$$
clearly satisfies
$$\lim_{x\to 0}f_\alpha(x)=0,\quad\lim_{x\to\infty}f_\alpha(x)=\infty,\quad f_\alpha'(x)=\cosh x(x+(1-\alpha)\tanh x).$$
By \eqref{condbarsigma}, we conclude that there is a unique $\rho_*>0$ with $f_\alpha(\rho_*)=0$. Now the triple $(\sigma_*,p_*,\rho_*)$ constitutes the unique flat stationary solution of \eqref{problem}.

The condition \eqref{condbarsigma} is also reasonable concerning the long-time behavior of our model: Let
$$\text{Vol}(\Omega_\rho(t))=\int_0^1\rho(t,x)\d x$$
denote the tumor volume. Using the condition $p_y(x,0)=0$ and the periodic boundary conditions for $p_x$, an application of the Gauss-Green Theorem shows that
$$\frac{\dw}{\dw t}\text{Vol}(\Omega_\rho(t))=\int_0^1\rho_t(t,x)\d x=-\int_{\partial\Omega_\rho(t)}\frac{\partial p}{\partial\nu}\d x=-\int_{\Omega_\rho(t)}\Delta p\,\d(x,y).$$
Applying the maximum principle, we conclude that
\beq\frac{\dw}{\dw t}\text{Vol}(\Omega_\rho(t))=
\mu\int_{\Omega_\rho(t)}(\sigma-\tilde\sigma)\d(x,y)\leq\mu(\max\{\bar\sigma_1,\bar\sigma_2\}-\tilde\sigma)\text{Vol}(\Omega_\rho(t)).\label{estimate}\eeq
Now \eqref{condbarsigma} guarantees that the right-hand side of \eqref{estimate} is positive; otherwise we would have $\text{Vol}(\Omega_\rho(t))\to 0$, as $t\to\infty$, meaning that the tumor will eventually vanish.

We have proved the following theorem.
\thm Assume that the condition \eqref{condbarsigma} is satisfied. Then the problem \eqref{problem} has a unique flat stationary solution $(\sigma_*,p_*,\rho_*)$ which is determined by the formulas \eqref{sigmastar}, \eqref{pstar} and \eqref{rhostar}.
\endthm\rm
\section{The linearization and asymptotic stability}\label{sec_stabil}
After studying local solvability and regularity of fully nonlinear equations, the second step is to consider the asymptotic behavior and in particular stability of the stationary solutions. Recall that a stationary solution $\bar u$ of an autonomous problem $u'(t)=F(u(t))$, $t>0$, is called stable if for each $\eps>0$ there is $\delta>0$ such that for $\Vert u_0-\bar u\Vert<\delta$, we have that $\Vert u(t)-\bar u\Vert<\eps$ for any $t>0$, and the solution $u=u(t;u_0)$ exists for all $t>0$ (denoted as $\tau(u_0)=\infty$). The stationary solution $\bar u$ is called asymptotically stable if it is stable and in addition $\Vert u(t)-\bar u\Vert\to 0$ as $t\to\infty$, uniformly for $u_0$ in a neighborhood of $\bar u$. It is said to be unstable if it is not stable, \cite{L95}.

In this section, it is our aim to prove that the stationary point $(\sigma_*,p_*,\rho_*)$ obtained in the previous section is asymptotically stable. Therefore, we consider the linearization of \eqref{problem2}: we plug the ansatz
$$
\begin{pmatrix}
  \tau \\
  q \\
  \rho \\
\end{pmatrix}
=
\begin{pmatrix}
  \sigma_*(y'\rho_*) \\
  p_*(y'\rho_*) \\
  \rho_* \\
\end{pmatrix}
+\eps
\begin{pmatrix}
  \Sigma(t,x',y') \\
  P(t,x',y') \\
  r(t,x')
\end{pmatrix}
$$
for $\eps>0$ small and with the new unknowns $(\Sigma,P,r)$ into \eqref{problem2} and compute the derivative with respect to $\eps$ at $\eps=0$. This yields
\beq
\left\{
\begin{array}{rclc}
  \Sigma_{x'x'}+\frac{1}{\rho_*^2}\Sigma_{y'y'}     & = & b(\sigma_*)r+\Sigma & \text{in }\Omega\times(0,T),\\
  P_{x'x'}+\frac{1}{\rho_*^2}P_{y'y'}               & = & b(p_*)r-\mu\Sigma   & \text{in }\Omega\times(0,T),\\
  r_t+\frac{1}{\rho_*}\left.P_{y'}\right|_{y'=1}    & = & 0                   & \text{on }\S\times(0,T),\\
  \Sigma|_{y'=1}                                    & = & 0                   & \text{on }\S\times(0,T),\\
  P|_{y'=1}                                         & = & -\gamma r_{x'x'}    & \text{on }\S\times(0,T),\\
  \Sigma|_{y'=0}                                    & = & 0                   & \text{on }\S\times(0,T),\\
  P_{y'}|_{y'=0}                                    & = & 0                   & \text{on }\S\times(0,T),\\
  r                                                 & = & r_0,                & t=0,
\end{array}
\right.
\label{problemlin}\eeq
where
$$b(v)r=\frac{2r}{\rho_*}v''(y'\rho_*)+r_{x'x'}y'v'(y'\rho_*),\quad v\in C^2[0,1].$$
For a given function $r\in h^{4+\alpha}(\S)$, we solve the boundary value problem for $\Sigma$ in \eqref{problemlin} and obtain a unique solution $\Sigma\in h^{4+\alpha}(\overline\Omega)$, which is periodic in $x$. Substituting $\Sigma$ into the second line we get a linear problem for $P$. Solving the equation for $P$ with the corresponding boundary conditions, we obtain a function $P\in h^{2+\alpha}(\overline\Omega)$ which is also periodic in $x$; observe that $\gamma r_{x'x'}\in h^{2+\alpha}(\S)$. Let us now introduce an operator $A\in\L(h^{4+\alpha}(\S),h^{1+\alpha}(\S))$ by setting
\beq(Ar)(x')=\frac{1}{\rho_*}\frac{\partial P}{\partial y'}(x',1),\quad x'\in\S.\label{opA}\eeq
Let
$$\Psi(\rho):=\Phi(\rho)\rho-F(\rho),\quad\rho\in h^{4+\alpha}(\S).$$
By \eqref{PhiF}, we have $\Psi\in C^{\infty}(h_+^{4+\alpha}(\S),h^{1+\alpha}(\S))$ and \eqref{abstrevoleqn} shows that $\rho_t=-\Psi(\rho)$ for any solution $\rho$ of \eqref{problem}. With
$$-D\Psi(\rho_*)r=-\left.\frac{\dw}{\dw\eps}\Psi(\rho_*+\eps r)\right|_{\eps=0}=\left.\frac{\dw}{\dw\eps}(\rho_*+\eps r)_t\right|_{\eps=0}=r_t$$
and the third equation in \eqref{problemlin} we see that $A=D\Psi(\rho_*)$.

As in \cite{CuiE09}, we can conclude the following proposition.
\prop\label{thm_Asemigroup} We have $A\in\L(h^{4+\alpha}(\S),h^{1+\alpha}(\S))$, and $-A$, considered as an unbounded operator in $h^{1+\alpha}(\S)$, generates a strongly continuous analytic semigroup.
\endprop\rm
It is our goal to represent the operator $A$ introduced in \eqref{opA} as a Fourier multiplication operator. To do so, we use that the functions $r\in h^{4+\alpha}(\S)$, $\Sigma\in h^{4+\alpha}(\overline\Omega)$ and $P\in h^{2+\alpha}(\overline\Omega)$ are periodic with respect to the variable $x'$.
We thus have the expansions
\bea
  \begin{pmatrix}
    r(x') \\
    \Sigma(x',y') \\
    P(x',y') \\
  \end{pmatrix}
=
  \begin{pmatrix}
    a_0 \\
    A_0(y') \\
    M_0(y') \\
  \end{pmatrix}
+
\sum_{k=1}^{\infty}
  \begin{pmatrix}
    a_k & b_k\\
    A_k(y') & B_k(y') \\
    M_k(y') & N_k(y') \\
  \end{pmatrix}
\cdot
  \begin{pmatrix}
  \cos (kx') \\
  \sin (kx') \\
  \end{pmatrix}.
\nonumber
\eea
\prop\label{propopA} The operator $A$ is a Fourier multiplication operator, i.e., given $r\in C^{\infty}(\S)$ with the Fourier expansion
$$r(x')=a_0+\sum_{k=1}^{\infty}(a_k\cos(kx')+b_k\sin(kx')),$$
the image $Ar$ has the Fourier expansion
$$(Ar)(x')=\lambda_0a_0+\sum_{k=1}^{\infty}\lambda_k\left(a_k\cos(kx')+b_k\sin(kx')\right)$$
with
\bea
\lambda_k & \!\!\!=\!\!\! & \frac{\mu (\bar\sigma_2\cosh\rho_*-\bar\sigma_1)\sqrt{1+k^2}}{\sinh\rho_*\sinh(\rho_*\sqrt{1+k^2})\cosh(\rho_*k)}\left(\cosh(\rho_*\sqrt{1+k^2})\cosh(\rho_*k)-1\right)
\nonumber\\
&&\quad+\left(\gamma k^2-\mu\tilde\sigma\rho_*-\mu\frac{\bar\sigma_2-\bar\sigma_1\cosh\rho_*}{\sinh\rho_*}\right)
k\tanh(\rho_*k)+\mu(\tilde\sigma-\bar\sigma_2).
\label{lambdak}
\eea
Moreover, the spectrum $\sigma(A)$ of $A$ is
$$\sigma(A)=\{\lambda_k;\;k\in\N_0\}$$
and $\sigma(A)\subset\R_+$ for sufficiently large $\gamma$.
\endprop\rm
\proof To simplify notation, we will write
$$c_1=\frac{\bar\sigma_2-\bar\sigma_1\cosh\rho_*}{\sinh\rho_*},\quad c_2=\bar\sigma_1,\quad c_3=c_1\cosh\rho_*+c_2\sinh\rho_*$$
in the sequel. First, we solve the boundary value problems
$$
\left\{
\begin{array}{rcl}
  -k^2A_k(y')+\frac{1}{\rho_*^2}A_k''(y') & = & A_k(y')+a_kf_k(y'), \\
  A_k(0) & = & 0, \\
  A_k(1) & = & 0,
\end{array}
\right.
$$
for $k\in\N_0$, and
$$
\left\{
\begin{array}{rcl}
  -k^2B_k(y')+\frac{1}{\rho_*^2}B_k''(y') & = & B_k(y')+b_kf_k(y'), \\
  B_k(0) & = & 0, \\
  B_k(1) & = & 0,
\end{array}
\right.
$$
for $k\in\N$, with
$$f_k(y')=\frac{2}{\rho_*}(c_1\sinh(y'\rho_*)+c_2\cosh(y'\rho_*))-k^2y'(c_1\cosh(y'\rho_*)+c_2\sinh(y'\rho_*)).$$
A lengthy and somewhat tedious computation shows that
\bea
A_k(y') & = & a_kc_1\left(y'\cosh(y'\rho_*)-\frac{\sinh(y'\rho_*\sqrt{1+k^2})}{\sinh(\rho_*\sqrt{1+k^2})}\cosh\rho_*\right) \nonumber\\
&& + a_kc_2\left(y'\sinh(y'\rho_*)-\frac{\sinh(y'\rho_*\sqrt{1+k^2})}{\sinh(\rho_*\sqrt{1+k^2})}\sinh\rho_*\right), \nonumber
\eea
and $B_k(y')$ is obtained from $A_k(y')$ by exchanging $a_k$ with $b_k$. The next task is to solve the boundary value problems
$$
\left\{
\begin{array}{rcl}
  -k^2M_k(y')+\frac{1}{\rho_*^2}M_k''(y') & = & -\mu A_k(y')+a_kg_k(y'), \\
  M_k'(0) & = & 0, \\
  M_k(1)  & = & \gamma k^2a_k,
\end{array}
\right.
$$
for $k\in\N_0$, and
$$
\left\{
\begin{array}{rcl}
  -k^2N_k(y')+\frac{1}{\rho_*^2}N_k''(y') & = & -\mu B_k(y')+b_kg_k(y'), \\
  N_k'(0) & = & 0, \\
  N_k(1)  & = & \gamma k^2b_k,
\end{array}
\right.
$$
for $k\in\N$, with
\bea
g_k(y') & = & \frac{2}{\rho_*}\mu(\tilde\sigma - c_1\sinh(y'\rho_*)-c_2\cosh(y'\rho_*))\nonumber\\
&&-k^2\mu (c_1y'+\tilde\sigma\rho_*(y')^2-c_2y'\sinh(y'\rho_*)-c_1y'\cosh(y'\rho_*)).
\nonumber
\eea
Again, it is straightforward to derive the solutions
\bea
M_0(y') & = & -\frac{a_0\mu c_3}{\sinh\rho_*}(y'\rho_*-\rho_*-\sinh(y'\rho_*))-a_0\mu\tilde\sigma\rho_*-a_0\mu c_1\nonumber\\
& & +a_0\mu c_1y'+a_0\mu\tilde\sigma\rho_*(y')^2-a_0\mu(c_1y'\cosh(y'\rho_*)+c_2y'\sinh(y'\rho_*))\nonumber
\eea
and, for $k\neq 0$,
\bea
M_k(y') & = & -\frac{a_k\mu c_3\sqrt{1+k^2}}{k\sinh(\rho_*\sqrt{1+k^2})}\big(\sinh(y'\rho_*k)-\tanh(\rho_*k)\cosh(y'\rho_*k)\big)
\nonumber\\
&&+a_k\mu c_3\frac{\sinh(y'\rho_*\sqrt{1+k^2})}{\sinh(\rho_*\sqrt{1+k^2})}+a_k(\gamma k^2-\mu\tilde\sigma\rho_*-\mu c_1)\frac{\cosh(y'\rho_*k)}{\cosh(\rho_*k)}\nonumber\\
&&+a_k\mu c_1y'+a_k\mu\tilde\sigma\rho_*(y')^2-a_k\mu (c_1y'\cosh(y'\rho_*)+c_2y'\sinh(y'\rho_*));\nonumber
\eea
the $N_k$ are obtained from $M_k$ by replacing $a_k$ with $b_k$. By the definition of $A$, we see that $A$ is a Fourier multiplication operator with
$$\lambda_k=\frac{1}{a_k\rho_*}M_k'(1).$$
Since $h^{4+\alpha}(\S)$ is compactly embedded into $h^{1+\alpha}(\S)$ and the resolvent set $\rho(A)$ is non-empty, the resolvent operator $(A-\lambda)^{-1}$ is compact for any $\lambda\notin\sigma(A)$. Hence the spectrum $\sigma(A)$ of $A$ consists entirely of eigenvalues.
Our explicit computations show that, for all $k\in\N_0$,
\bea
\lambda_k & \!\!=\!\! & -\frac{\mu c_3\sqrt{1+k^2}}{\sinh(\rho_*\sqrt{1+k^2})\cosh(\rho_*k)}+\frac{\mu c_3\sqrt{1+k^2}}{\tanh(\rho_*\sqrt{1+k^2})}\nonumber\\
&&\quad +(\gamma k^2-\mu\tilde\sigma\rho_*-\mu c_1)k\tanh(\rho_*k)\nonumber\\
&&\quad +2\mu\tilde\sigma+\frac{\mu}{\rho_*}(c_1-c_3)-\mu(c_1\sinh\rho_*+c_2\cosh\rho_*),\nonumber
\eea
and, by our definitions, $\lambda_k$ is as specified in \eqref{lambdak}. Concerning positivity, we first concentrate our attention on $\lambda_0$. Using once again \autoref{rhostar}, we obtain
$$\lambda_0=\mu\tilde\sigma\left(1-\frac{\rho_*}{\sinh\rho_*}\right),$$
and written in this form, $\lambda_0$ is clearly positive. For $k\neq 0$, fix some $\gamma_0>0$ and let $\lambda_k(\gamma_0)$ be as in \eqref{lambdak}. Since
$$\frac{\lambda_k(\gamma_0)}{k^3\tanh(\rho_*k)}=\gamma_0+n_k,\quad n_k\to 0\,\text{ for }\,k\to\infty,$$
we see that there exists $k_0\in\N$ such that $\lambda_k(\gamma_0)>0$ for all $k>k_0$. Since $\lambda_k(\gamma)\geq\lambda_k(\gamma_0)$ for $\gamma\geq\gamma_0$, it follows that $\lambda_k(\gamma)>0$ for all $k>k_0$ and $\gamma\geq\gamma_0$. It is obvious that we can choose $\gamma\geq\gamma_0$ so large that also $\lambda_1(\gamma),\ldots,\lambda_{k_0}(\gamma)>0$. Let $\Lambda\neq\emptyset$ be the collection of all $\gamma>0$ such that $\lambda_k(\gamma)>0$ for all $k\in\N$. We have shown that $\sigma(A)\subset\R_+$ for any $\gamma\in\Lambda$ and this achieves the proof.
\endproof
\begin{figure}[H]
\begin{center}
\includegraphics[width=9cm]{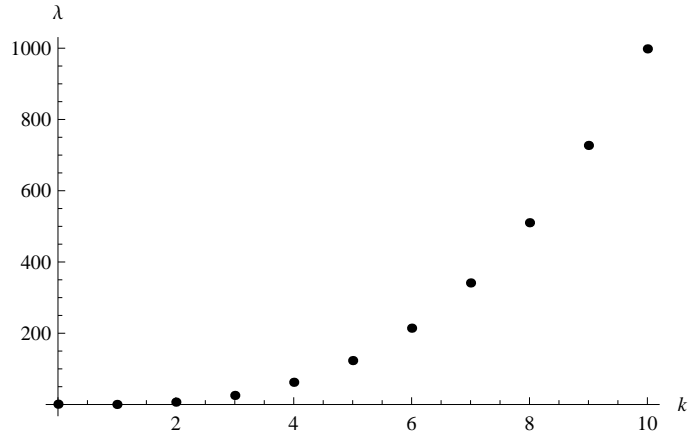}
\caption{The spectrum of the operator $A$, for $\mu=\tilde\sigma=\gamma=1$, $\bar\sigma_1=2$, $\bar\sigma_2=3$.}\label{fig1}
\end{center}
\end{figure}
In fact, not the operator $A$ but $-A$ will determine the stability properties of the problem \eqref{problem}. We will employ the following stability theorem.
\thm[see \cite{L95}]\label{thm_stab} Let $X$ be a Banach space and let $T\colon D(T)\subset X\to X$ be a linear sectorial operator such that the graph norm of $T$ is equivalent to the norm of $D$. Assume furthermore that
$$
\sup\{\text{\rm Re }\lambda,\,\lambda\in\sigma(T)\}=-\omega_0<0.
$$
Let $\mathcal O$ be a neighborhood of zero in $D$ and let $G\colon\mathcal O\to X$ be a $C^1$ function with locally Lipschitz continuous derivative such that $G(0)=0$ and $G'(0)=0$. Fix $\omega\in[0,\omega_0)$. Then there exist $r>0$ and $M>0$ such that for each $u_0\in B(0,r)\subset D$ the solution $u(t)$ to the problem
$$u'(t)=Tu(t)+G(u(t)),\,t>0,\quad u(0)=u_0,$$
satisfies $\tau(u_0)=\infty$ and
$$\Vert u(t)\Vert_D+\Vert u'(t)\Vert_X\leq Me^{-\omega t}\Vert u_0\Vert_D,\quad t\geq 0.$$
\endthm\rm
As for the model in \cite{CuiE09}, we show that the steady state $\rho_*$ of
\beq\rho_t+\Psi(\rho)=0,\quad\rho(0)=\rho_0,\label{IVPrho}\eeq
is asymptotically stable under small perturbations belonging to $h_+^{4+\alpha}(\S)$. By definition, $\Psi(\rho)=\Phi(\rho)\rho-F(\rho)$, for $\rho\in h_+^{4+\alpha}(\S)$. Setting $G(r):=\Psi(r+\rho_*)-Ar$, we have that
$$G\in C^{\infty}(h_+^{4+\alpha}(\S),h^{1+\alpha}(\S))$$
and observe that
$$G(0)=\Psi(\rho_*)=0,\quad DG(0)=D\Psi(\rho_*)-A=0.$$
Theorem~\ref{thm_stab} implies that the solution $r$ to
$$r_t=-Ar-G(r)$$
is asymptotically stable: There are constants $\omega$, $\eps$, $K$ such that if $r_0\in h_+^{4+\alpha}(\S)$ satisfies $\Vert r_0\Vert_{C^{4+\alpha}}<\eps$ then the solution $r$ exists globally and satisfies
$$\Vert r(t)\Vert_{C^{4+\alpha}}\leq K\exp(-\omega t)\Vert r_0\Vert_{C^{4+\alpha}},\quad t\geq 0.$$
Letting $r(t)=\rho(t)-\rho_*$ for $t\geq 0$, this in turn shows that if $\Vert\rho_0-\rho_*\Vert_{C^{4+\alpha}}<\eps$, then the solution to \eqref{IVPrho} exists globally and satisfies
\beq\Vert \rho(t)-\rho_*\Vert_{C^{4+\alpha}}\leq K\exp(-\omega t)\Vert \rho_0-\rho_*\Vert_{C^{4+\alpha}},\quad t\geq 0.\label{estimaterho}\eeq
Next let us consider the dynamical behavior of $\sigma$ and $p$. By definition, we have
$$\sigma_*=\mathcal R(\rho_*)(\bar\sigma_1,\bar\sigma_2)\quad\text{and}\quad\sigma(t)=\mathcal R(\rho(t))(\bar\sigma_1,\bar\sigma_2).$$
By the mean value theorem, there exists a constant $C$ such that
$$\Vert\sigma(t)-\sigma_*\Vert_{C^{4+\alpha}}=\Vert[\mathcal R(\rho(t))-\mathcal R(\rho_*)](\bar\sigma_1,\bar\sigma_2)\Vert_{C^{4+\alpha}}\leq C\Vert\rho(t)-\rho_*\Vert_{C^{4+\alpha}}$$
for any $t\geq 0$. Combining this with estimate \eqref{estimaterho}, we get
$$\Vert\sigma(t)-\sigma_*\Vert_{C^{4+\alpha}}\leq K\exp(-\omega t).$$
A corresponding estimate for $p$ can be obtained similarly. We have proved the following theorem.
\thm Let $\bar\sigma_1,\bar\sigma_2>\tilde\sigma>0$, $\mu>0$ and $\gamma\in\Lambda$ be given. Then the flat stationary solution defined by \eqref{sigmastar}, \eqref{pstar} and \eqref{rhostar} is asymptotically stable: There are positive constants $\omega$, $K$ and $\eps$ such that if $\Vert \rho_{0}-\rho_*\Vert_{C^{4+\alpha}} < \eps$ then
$$\Vert\sigma(t)-\sigma_*\Vert_{C^{4+\alpha}}+\Vert p(t)-p_*\Vert_{C^{2+\alpha}}+\Vert\rho(t)-\rho_*\Vert_{C^{4+\alpha}}\leq K\exp(-\omega t),$$
for any $t\geq 0$.
\endthm\rm

\end{document}